\documentclass{article}\begin{document}\centerline{\bf\ Notes On a Continued Fraction of Ramanujan}\vskip .10in

\centerline{\bf Nikos Bagis}

\centerline{Department of Informatics}

\centerline{Aristotele University of Thessaloniki Greece}

\centerline{nikosbagis@hotmail.gr}

\begin{quote}
\begin{abstract}
We study the properties of a general continued fraction of Ramanujan. In some certain cases we evaluate it completely.  
\end{abstract}

\bf keywords \rm{Continued Fractions; Ramanujan;}

\end{quote}

\section{Introduction}
\label{intro}
Let
\begin{equation}
\left(a;q\right)_k=\prod^{k-1}_{n=0}(1-aq^n)
\end{equation}
Then we define
\begin{equation}
f(-q)=(q;q)_\infty 
\end{equation}
and
\begin{equation}
\Phi(-q)=(-q;q)_\infty
\end{equation}
Also let
\begin{equation}
K(x)=\int^{\pi/2}_{0} \frac{1}{\sqrt{1-x^2\sin^2(t)}}dt
\end{equation}
be the elliptic integral of the first kind\\
The function $k_r$ is defined from the equation  
\begin{equation}
\frac{K(k'_r)}{K(k_r)}=\sqrt{r}
\end{equation}
where $r$ is positive , $q=e^{-\pi \sqrt{r}}$ and $k'=\sqrt{1-k^2}$. Note also that whenever $r$ is positive rational, the $k$ are algebraic numbers.
 
In Berndt's book: Ramanujan's Notebook Part III,  ([B3] pg.21), one can find the following expansion
\[
\]
\textbf{Theorem.}\\
Suppose that either $q, a$ and $b$ are complex numbers with $\left|q\right|<1$, or $q,a$, and $b$ are complex numbers with $a=bq^m$ for some integer $m$. Then
$$U=U(a,b;q)=\frac{(-a;q)_{\infty}(b;q)_{\infty}-(a;q)_{\infty}(-b;q)_{\infty}}{(-a;q)_{\infty}(b;q)_{\infty}+(a;q)_{\infty}(-b;q)_{\infty}}=$$  
\begin{equation}
\frac{a-b}{1-q+}\frac{(a-bq)(aq-b)}{1-q^3+}\frac{q(a-bq^2)(aq^2-b)}{1-q^5+}\frac{q^2(a-bq^3)(aq^3-b)}{1-q^7+}\ldots
\end{equation}
\[
\]
Suppose now 
\begin{equation}
X=\frac{(-a;q)_{\infty}(b;q)_{\infty}}{(a;q)_{\infty}(-b;q)_{\infty}}
\end{equation}
Then holds
\begin{equation}
\frac{X-1}{X+1}=U
\end{equation}
\section{Propositions}
\textbf{Proposition 1.}\\
Set 
\begin{equation}
\phi(q)=\sum^{\infty}_{n=-\infty}q^{n^2}
\end{equation}
then
\begin{equation}
\frac{\phi(q)-1}{\phi(q)+1}=\frac{q}{1+q+}\frac{-q^3}{1+q^3+}\frac{-q^5}{1+q^5+}\frac{-q^7}{1+q^7+}\ldots
\end{equation}
\textbf{Proof.}\\
Take $q\rightarrow q^2$ in (6) and then set $a\rightarrow q$ and $b\rightarrow q^2$.
\[
\]
\textbf{Proposition 2.}
\begin{equation}
\frac{\Phi(-q)-f(-q)}{\Phi(-q)+f(-q)}=\frac{q}{1-q+}\frac{q^3}{1-q^3+}\frac{q^5}{1-q^5+}\frac{q^7}{1-q^7+}\ldots
\end{equation}
\textbf{Proof.}\\
Set $b=0$ in (6) and then $a=q$.
\[
\]
\textbf{Proposition 3.}
\begin{equation}
\frac{\Phi(-q)-f(-q)}{\Phi(-q)+f(-q)}=-\frac{\phi(-q)-1}{\phi(-q)+1}
\end{equation}
\textbf{Proof.}\\
It follows from Propositions 1, 2
\[
\]
\textbf{Proposition 4.}
\begin{equation}
\sum^{\infty}_{n=0}\frac{q^n}{1-a^2q^{2n}}=\frac{1}{1-q+}\frac{-a^2(1-q)^2}{1-q^3+}\frac{-qa^2(1-q^2)^2}{1-q^5+}\frac{-q^2a^2(1-q^3)^2}{1-q^7+}\ldots
\end{equation}
\textbf{Proof.}\\
Divide relation (6) by $a-b$ and then take the limit $b\rightarrow a$.
\[
\]
\textbf{Proposition 5.}
\begin{equation}
\frac{K(k_r)}{2\pi}+\frac{1}{4}=\frac{1}{1-q+}\frac{(1-q)^2}{1-q^3+}\frac{q(1-q^2)^2}{1-q^5+}\frac{q^2(1-q^3)^2}{1-q^7+}\ldots
\end{equation}
\textbf{Proof.}\\
Set in (13) $a=i$, and $q=e^{-\pi \sqrt{r}}$.
\[
\]
Now set 
\begin{equation}
u(a,q)=\frac{2a}{1-q+}\frac{a^2(1+q)^2}{1-q^3+}\frac{a^2q(1+q^2)^2}{1-q^5+}\frac{a^2q^2(1+q^3)^2}{1-q^7+}\ldots
\end{equation}
and
\begin{equation}
P=\left(\frac{(-a;q)_{\infty}}{(a;q)_{\infty}}\right)^2
\end{equation}
Then
\begin{equation}
\frac{P-1}{P+1}=u(a,q)
\end{equation}
or
\[
\]
\textbf{Proposition 6.}
\begin{equation}
\left(\frac{(-a;q)_{\infty}}{(a;q)_{\infty}}\right)^2=-1+\frac{2}{1-}\frac{2a}{1-q+}\frac{a^2(1+q)^2}{1-q^3+}\frac{a^2q(1+q^2)^2}{1-q^5+}\frac{a^2q^2(1+q^3)^2}{1-q^7+}\ldots
\end{equation}
\[
\]
\textbf{Proposition 7.}
\begin{equation} 
4\sum^{\infty}_{n=0}\frac{a^{2n+1}}{(2n+1)(1-q^{2n+1})}=\log \left(-1+\frac{2}{1-u(a,q)}\right)
\end{equation}
\textbf{Proof.}\\
Take the logarithms in both sides of (18) and expand in Taylor series. Then rearange the double sum to get easily the desired result.
\[
\]   
Here we must mention that holds the more general formula
\begin{equation}
2\sum^{\infty}_{n=0}\frac{a^{2n+1}-b^{2n+1}}{(2n+1)(1-q^{2n+1})}=\log\left(-1+\frac{2}{1-U(a,b;q)}\right)
\end{equation}
Thus
\begin{equation}
\left(-1+\frac{2}{1-U(a,b;q)}\right)=\frac{\left(-1+\frac{2}{1-u(a,q)}\right)}{\left(-1+\frac{2}{1-u(b,q)}\right)}
\end{equation}
and for to study $U$ we have to study only $u$.   
In some cases the $u$ fraction can calculated in terms of elliptic functions. For example:
\begin{equation}
-1+\frac{2}{1-u(q,q)}=\frac{\pi}{2k'_rK(k_r)}
\end{equation}
In general holds
$$4\sum^{\infty}_{n=0}\frac{q^{\nu(2n+1)}}{(2n+1)(1-q^{2n+1})}=-4\sum^{\nu-1}_{j=1}\textrm{arctanh}(q^j)-\log\left(\frac{2k'_rK(k_r)}{\pi}\right)$$
from which we lead to the following:
\[
\]
\textbf{Proposition 8.}\\
Let $\nu$ be positive integer, then
$$
-1+\frac{2}{1-u(q^{\nu},q)}=-1+\frac{2}{1-}\frac{2q^{\nu}}{1-q+}\frac{q^{2\nu}(1+q)^2}{1-q^3+}\frac{q^{2\nu+1}(1+q^2)^2}{1-q^5+}\frac{q^{2\nu+2}(1+q^3)^2}{1-q^7+}\ldots=$$
\begin{equation}
=\frac{\pi}{2k'_rK(k_r)}\exp\left(-4\sum^{\nu-1}_{j=1}\textrm{arctanh}(q^j)\right)
\end{equation} 
\[
\]
\textbf{Proposition 9.}\\
Let $\nu_1$, $\nu_2$ be positive integers, then
\begin{equation}
-1+\frac{2}{1-U(q^{\nu_1},q^{\nu_2},q)}= \exp\left(-4\left(\sum^{\nu_1-1}_{j_1=1}\textrm{arctanh}(q^{j_1})-\sum^{\nu_2-1}_{j_2=1}\textrm{arctanh}(q^{j_2})\right)\right)
\end{equation}
\textbf{Proof.}\\ The proof follows easily from (21) and (24).
\[
\]
\textbf{Note.}
One can find many useful results in pages stored on the Web one is:
\[
\] http://pi.physik.uni-bonn.de/~dieckman/InfProd/InfProd.html 
\[
\] 
Another formula related with $u$ continued fraction is when $q=e^{-\pi\sqrt{r}}$   
$$
-1+\frac{2}{1-u(q^{\nu+1/2},q)}=\exp\left(-4\sum^{\infty}_{n=0}\frac{q^{(2n+1)(\nu+1/2)}}{(2n+1)(1-q^{2n+1})}\right)
$$
\begin{equation}
=
\exp\left(-4\sum^{\nu-1}_{j=0}\textrm{arctanh}(q^{j+1/2})+\textrm{arctanh}(k_r)\right)
\end{equation}
Hence also
$$
k_r=\tanh\left(4\sum^{\nu-1}_{j=0}\textrm{arctanh}(q^{j+1/2})+\log\left(-1+\frac{2}{1-u(q^{\nu+1/2},q)}\right)\right)
$$
For every $\nu$ positive integer.\\
Hence we obtain a continued fraction for $k_r$
\begin{equation}
\frac{k'_r}{1-k_r}=-1+\frac{2}{1-u(q^{1/2},q)}
\end{equation} 
\[
\]
Inspired from the above relations and Propositions we have
\[
\]
\textbf{Theorem}(Unproved)\\ 
If $c$ is positive real and $\nu_1, \nu_2$ positive integers then:
$$-1+\frac{2}{1-U(q^{\nu_1+c},q^{\nu_2+c},q)}\stackrel{?}{=}$$ 
\begin{equation}
=\exp\left(-4\left(\sum^{\nu_1-1}_{j_1=1}\textrm{arctanh}(q^{j_1+c})-\sum^{\nu_2-1}_{j_2=1}\textrm{arctanh}(q^{j_2+c})\right)\right)
\end{equation}
or better\\
If $U=U(a,b,q)$, where $q=e^{-\pi\sqrt{r}}$ and 
\begin{equation}
c=\left\{\frac{-\log(a)}{\pi\sqrt{r}}\right\}=\left\{\frac{-\log(b)}{\pi\sqrt{r}}\right\}
\end{equation}
then 
$$-1+\frac{2}{1-U(a,b,q)}=$$
$$=\exp\left(-4\left(\sum^{\left[\frac{-\log(a)}{\pi\sqrt{r}}\right]-1}_{j_1=1}\textrm{arctanh}(q^{j_1+c})-\sum^{\left[\frac{-\log(b)}{\pi\sqrt{r}}\right]-1}_{j_2=1}\textrm{arctanh}(q^{j_2+c})\right)\right)$$
where $\left\{x\right\}$ is the fractional part of $x$ and $\left[x\right]$ is the largest integer that not exiding $x$.\\   
ii) Observe that for $\nu_1=\nu_2=\nu$ $$-1+\frac{2}{1-U(q^{\nu+c},q^{\nu+c},q)}=1$$
iii) Also we observe that holds and
\begin{equation}
\left(-1+\frac{2}{1-U(a,-b;q)}\right)=\left(-1+\frac{2}{1-u(a,q)}\right)\left(-1+\frac{2}{1-u(b,q)}\right)
\end{equation}
This relation is similarly to (21). We also get the following unproved
\[
\]
\textbf{Proposition 10.}(Unproved)\\ 
Let $w\in \textrm{Im}(\bf C\rm)$, then
\begin{equation}
\left|-1+\frac{2}{1-U(q^{\nu_1+c},-w  q^{\nu_2+c},q)}\right|\stackrel{?}{=} \left(-1+\frac{2}{1-u(q^{\nu_1+c},q)}\right)
\end{equation}
\[
\]
Seting $c=0$ in (30) and using Proposition 9 we get
\[
\] 
\textbf{Proposition 11.}\\ When $w, z\in \textrm{Im}(\bf C\rm)$ and $q=e^{-\pi\sqrt{r}}$, $r>0$, then\\
(i)
$$
\left|-1+\frac{2}{1-U(q^{\nu_1},-w  q^{\nu_2},q)}\right|=
\frac{\pi}{2k'_rK(k_r)}\exp\left(-4\sum^{\nu_1-1}_{j=1}\textrm{arctanh}(q^j)\right)
$$
(ii)
$$
\left|-1+\frac{2}{1-U(-z q^{\nu_1+c},-w q^{\nu_2+c},q)}\right|=1
$$

\newpage

\centerline{\bf References}\vskip .2in

[1]:M.Abramowitz and I.A.Stegun, Handbook of Mathematical Functions. Dover Publications

[2]:B.C.Berndt, Ramanujan`s Notebooks Part I. Springer Verlang, New York (1985)

[3]:B.C.Berndt, Ramanujan`s Notebooks Part II. Springer Verlang, New York (1989)

[4]:B.C.Berndt, Ramanujan`s Notebooks Part III. Springer Verlang, New York (1991)

[5]:L.Lorentzen and H.Waadeland, Continued Fractions with Applications. Elsevier Science Publishers B.V., North Holland (1992) 

[6]:E.T.Whittaker and G.N.Watson, A course on Modern Analysis. Cambridge U.P. (1927)

\end{document}